\documentclass[11pt]{amsart}
\usepackage{amssymb,amsmath,epsfig,graphics}
\theoremstyle{plain}
\newtheorem{thrm}{Theorem}[section]
\newtheorem{lemma}[thrm]{Lemma}
\newtheorem{prop}[thrm]{Proposition}
\newtheorem{cor}[thrm]{Corollary}

\newtheorem{dfn}[thrm]{Definition}

\numberwithin{equation}{section} \numberwithin{figure}{section}

\begin{document}
\newcommand{\SL}{\mathcal L^{1,p}(\Om)}
\newcommand{\Lp}{L^p(\Omega)}
\newcommand{\CO}{C^\infty_0(\Omega)}
\newcommand{\Rn}{\mathbb R^n}
\newcommand{\Rm}{\mathbb R^m}
\newcommand{\R}{\mathbb R}
\newcommand{\Om}{\Omega}
\newcommand{\Hn}{\mathbb H^n}
\newcommand{\HH}{\mathbb H^1}
\newcommand{\eps}{\epsilon}
\newcommand{\BVX}{BV_H(\Omega)}
\newcommand{\IO}{\int_\Omega}
\newcommand{\bG}{\boldsymbol{G}}
\newcommand{\bg}{\mathfrak g}
\newcommand{\p}{\partial}
\newcommand{\Xnu}{\overset{\rightarrow}{ X_\nu}}
\newcommand{\nuX}{\boldsymbol{\nu}_H}
\newcommand{\Up}{\boldsymbol{\mathcal Y}_H}
\newcommand{\n}{\boldsymbol \nu}
\newcommand{\sigmau}{\boldsymbol{\sigma}^u_H}
\newcommand{\nui}{\nu_{H,i}}
\newcommand{\nuj}{\nu_{H,j}}
\newcommand{\dej}{\delta_{H,j}}
\newcommand{\cx}{\boldsymbol{c}_\mathcal S}
\newcommand{\sx}{\sigma_H}
\newcommand{\lx}{\mathcal L_H}
\newcommand{\pb}{\overline p}
\newcommand{\qb}{\overline q}
\newcommand{\ob}{\overline \omega}
\newcommand{\nuu}{\boldsymbol \nu_{H,u}}
\newcommand{\nuv}{\boldsymbol \nu_{H,v}}
\newcommand{\Bl}{\Bigl|_{\lambda = 0}}
\newcommand{\mS}{\mathcal S}
\newcommand{\delh}{\Delta_H}
\newcommand{\delinf}{\Delta_{H,\infty}}
\newcommand{\nabh}{\nabla^H}
\newcommand{\delp}{\Delta_{H,p}}
\newcommand{\mO}{\mathcal O}
\newcommand{\delhs}{\Delta_{H,S}}
\newcommand{\lhs}{\hat{\Delta}_{H,S}}
\newcommand{\bN}{\boldsymbol{N}}
\newcommand{\bnu}{\boldsymbol \nu}
\newcommand{\la}{\lambda}
\newcommand{\nup}{\boldsymbol{\nu}_H^\perp}
\newcommand{\fv}{\mathcal V^{H}_I(\mS;\mathcal X)}
\newcommand{\sv}{\mathcal V^{H}_{II}(\mS;\mathcal X)}
\newcommand{\di}{\nabla_{i}^{H,\mS}}
\newcommand{\one}{\nabla_{1}^{H,\mS}}
\newcommand{\two}{\nabla_{2}^{H,\mS}}
\newcommand{\del}{\nabla^{H,\mS}}
\newcommand{\delXY}{\nabla^{H,\mS}_X Y}
\newcommand{\oX}{\overline X}
\newcommand{\oY}{\overline Y}
\newcommand{\ou}{\overline u}
\newcommand{\duno}{\nabla^{H,\mS}_1}
\newcommand{\ddue}{\nabla^{H,\mS}_2}

\begin{abstract} One of the main objectives
of this paper is to unravel a new interesting phenomenon of the
sub-Riemannian Bernstein problem with respect to its Euclidean
ancestor, with the purpose of also indicating a possible line of
attack toward its solution. We show that the global intrinsic graphs
\eqref{ce2i} are unstable critical points of the horizontal
perimeter. As a consequence of this fact, the study of the stability
acquires a central position in the problem itself.
\end{abstract}


\title[A notable family of entire
intrinsic minimal graphs, etc.] {A notable family of entire
intrinsic minimal graphs in the Heisenberg group which are not
perimeter minimizing}

\author{D. Danielli}
\address{Department of Mathematics\\Purdue University \\
West Lafayette, IN 47907} \email[Donatella
Danielli]{danielli@math.purdue.edu}
\thanks{First author supported in part by NSF CAREER Grant, DMS-0239771}

\author{N. Garofalo}
\address{Department of Mathematics\\Purdue University \\
West Lafayette, IN 47907} \email[Nicola
Garofalo]{garofalo@math.purdue.edu}
\thanks{Second author supported in part by NSF Grant DMS-0300477}

\author{D. M. Nhieu}
\address{Department of Mathematics\\Georgetown University \\
Washington, DC 20057-1233} \email[Duy-Minh
Nhieu]{nhieu@math.georgetown.edu}

%
%
 \subjclass{}
\date{\today}

\maketitle

\baselineskip 13.5pt



\section{\textbf{Introduction}}

\vskip 0.2in

The development of geometric measure theory in sub-Riemannian
spaces has received a strong impulse over the past decade, see
 \cite{Pa1}, \cite{Pa2}, \cite{CDG},
\cite{KR}, \cite{E1}, \cite{E2}, \cite{E3}, \cite{Gro}, \cite{GN},
\cite{Be}, \cite{DS}, \cite{DGN1}, \cite{AK1}, \cite{AK2},
\cite{CS1}, \cite{A}, \cite{FSS1}, \cite{Ma1}, \cite{FSS2},
\cite{Ma2}, \cite{CMS}, \cite{FSS3}, \cite{BRS}, \cite{DGN4},
\cite{DGN5}, \cite{LR}, \cite{LM}, \cite{FSS4}, \cite{Ma3},
\cite{CS2} \cite{P1}, \cite{P2}, \cite{GP1}, \cite{CG}, \cite{CHMY},
\cite{CH}, \cite{HP}, \cite{HP2}, \cite{RR}, \cite{BC}, \cite{Se},
\cite{Se2}, \cite{Mo}. In particular, the papers \cite{GP1},
\cite{CHMY} and \cite{CH} contain a detailed study of the Bernstein
problem in the first Heisenberg group $\HH$, and in more general CR
manifolds of real dimension three. Despite the progress made in
these latter papers, this problem presently still constitutes a
basic open question. One of the main objectives of this paper is to
unravel a new interesting phenomenon of the sub-Riemannian Bernstein
problem with respect to its Euclidean ancestor, with the purpose of
also indicating a possible line of attack toward its solution.

To provide the reader with some perspective, we recall that in the
Heisenberg group $\Hn$ a basic discovery of Franchi, Serapioni and
Serra Cassano is a structure theorem \`a la De Giorgi for sets of
locally finite horizontal perimeter \cite{FSS1} (see \cite{FSS3} for
an extension to Carnot groups of step $r=2$). To prove the latter
they show that the non-isotropic blow-up of such a set at a point of
its reduced boundary produces a vertical hyperplane
\begin{equation}\label{vpi}
P_\gamma\ =\ \{(x,y,t)\in \Hn \mid <a,x> + <b,y> = \gamma\}\ ,\
\quad\quad\quad a^2 + b^2 \not= 0\ ,
\end{equation}
(when $\gamma = 0$ these sets are also the maximal subgroups of
$\Hn$). Recalling that the characteristic locus of a hypersurface
$\mS\subset \Hn$, denoted henceforth by $\Sigma(\mS)$, is the
collection of all points $g\in \mS$ at which $T_g\mS = H_g\Hn$,
where $H\Hn$ denotes the horizontal bundle of $\Hn$, it is easy to
recognize that $\Sigma(P_\gamma) = \varnothing$ for any $\gamma \in
\R$. Therefore, in analogy with the Euclidean situation, the cited
blow-up result from \cite{FSS1} suggests the natural conjecture that
if $\mS\subset \Hn$ is a $C^2$ entire $H$-minimal graph over some
hyperplane, and if $\Sigma(\mS)= \varnothing$, then $\mS = P_\gamma$
for some $\gamma$. Here, following a perhaps unfortunate tradition
of the classical situation, $H$-minimal is intended in the sense
that $\mS$ is of class $C^2$, and the horizontal mean curvature,
defined in \eqref{HMCgen} below, vanishes identically as a
continuous function on $\mS$. However, for $\Hn$ the situation is
very different than in Euclidean space. In fact, it was proved in
\cite{GP1} that there exist non-planar entire $H$-minimal graphs
with empty characteristic locus, thus violating the above plausible
conjecture. For instance, the non-planar real analytic surfaces
\begin{equation}\label{ce2i}
\mS\ =\ \{(x,y,t)\in \HH\mid x = y (\alpha  t + \beta) , \alpha >0,
\beta \in \R\}\ ,
\end{equation}
are $H$-minimal, and they have empty characteristic locus. Moreover,
the surfaces \eqref{ce2i} are global intrinsic $X_1$-graphs in the
sense of Franchi, Serapioni and Serra Cassano, see \cite{FSS4} (we
stress that intrinsic graphs have empty characteristic locus by
definition). This means that there exists a globally defined
function $\phi:\R^2_{u,v} \to \R$ such that, in the coordinates
$(x,y,t)$, we can parameterize \eqref{ce2i} as follows
\begin{equation}\label{X1graph}
(x,y,t) = (0,u,v) \circ \phi(u,v) e_1 = (0,u,v)\ \circ
(\phi(u,v),0,0) = \left(\phi(u,v),u,v - \frac{u}{2}
\phi(u,v)\right)\ ,
\end{equation}
where $\circ$ indicates the non-Abelian group multiplication in
$\HH$, see \eqref{gl}. In fact, imposing the defining equation $x =
y(\alpha t + \beta)$ for $\mS$, with $\alpha
>0$, we obtain the function
\begin{equation}\label{phi1}
\phi(u,v)\ =\ \frac{2 u (\alpha v + \beta)}{2 + \alpha u^2}\ ,
\end{equation}
which describes $\mS$ as an entire $X_1$-graph. We also note that a
vertical plane $P_\gamma$ is a global $X_1$-graph if $a \not= 0$ (or
an $X_2$-graph if $b\not= 0$), with
\begin{equation}\label{vp}
\phi(u,v)\ =\ \frac{\gamma - b u}{a}\ ,\quad\text{if}\quad a\not= 0\
.
\end{equation}

Examples such as \eqref{ce2i} seem to cast a dim light on the
Bernstein problem in $\HH$. There is however a deeper aspect of the
problem which has gone unnoticed so far. What could be happening in
fact is that, due to the different nature of the relevant perimeter
functional, global intrinsic graphs such as \eqref{ce2i} are only
stationary for the horizontal perimeter, but not stable. In this
paper we examine this aspect in depth. More precisely, we recall
that, thanks to the convexity of the area functional
\[
\mathcal A(u)\ =\ \int_\Om \sqrt{1 +|Du|^2}\ dx\ ,\quad\quad\quad
\Om \subset \Rn\ ,
\] in the classical theory of minimal surfaces any critical point of
$\mathcal A$ is automatically stable, i.e., it is a local minimizer,
see e.g. \cite{CM}. By contrast, we show that the global intrinsic
graphs \eqref{ce2i} are unstable critical points of the horizontal
perimeter $P_H(\mS)$ defined in \eqref{sigmah} below. As a
consequence of this fact, the study of the stability acquires a
central position in the problem itself, and our results suggest
that, if properly understood from this new perspective, the
Bernstein property is still true in $\HH$. Besides their intrinsic
interest, we believe that the relevance of our results lies in the
method of proof, which is quite general and flexible, and has the
potential of being successfully applied to attack the sub-Riemannian
Bernstein problem. In Geometric Measure Theory there exist in
essence two (different, but equivalent) approaches to stability: the
former is based on the so-called method of calibrations, the latter
on second variation formulas. Our approach revolves around a general
second variation formula established in \cite{DGN3}, and on the
explicit construction of a continuum of directions along which the
intrinsic perimeter strictly decreases. We need to introduce a basic
definition.

\medskip

\begin{dfn}\label{D:minimali}
We say that a $C^2$ oriented $H$-minimal surface $\mS \subset \HH$,
with $\Sigma(\mS)=\varnothing$, is \emph{stationary} if it has
vanishing first variation of the $H$-perimeter, i.e., if
\[
\fv\ \overset{def}{=}\ \frac{d}{d\lambda}P_H(\mathcal
S^\lambda)\Bigl|_{\lambda = 0}\ =\ 0\ ,
\]
for any deformation $\mS \to \mS^\lambda = \mS + \lambda \mathcal
X$, where $\mathcal X \in C^2_0(\mS,\HH)$, with $\mathcal X
\not\equiv 0$. We say that a stationary $\mS$ is \emph{stable} if
the second variation is nonnegative, i.e.,
\[
\sv\ \overset{def}{=}\ \frac{d^2}{d\lambda^2}P_H(\mathcal
S^\lambda)\Bigl|_{\lambda = 0}\ \geq\ 0\ ,
\]
for any $\mathcal X \in C^2_0(\mS,\HH)$, with $\mathcal X \not\equiv
0$. If there exists such an $\mathcal X$ for which $\sv < 0$, then
we say that $\mS$ is \emph{unstable}.
\end{dfn}

\medskip

It has been proved in \cite{DGN3} that $\mS$ is stationary if and
only if $\mS$ is $H$-minimal, see Theorem \ref{T:variations} below.
In this paper we establish the following result.

\medskip

\begin{thrm}\label{T:permin}
For every $\alpha > 0$ and $\beta \in \R$, the $H$-minimal global
intrinsic $X_1$-graphs
\begin{equation}\label{min}
\mathcal S = \{(x,y,t)\in \HH \mid x = y (\alpha t + \beta) \}\ ,
\end{equation}
are unstable. More precisely, there exist $a\in C^\infty_0(\mS)$,
$a\not\equiv 0$, and $h\in C^\infty_0(\mS)$, $h\not\equiv 0$, such
that with either $\mathcal X = a X_1$, or $\mathcal X = h \nuX$, we
have $\sv < 0$.
\end{thrm}

\medskip

To explain the strategy behind Theorem \ref{T:permin} we mention
that in Section \ref{S:prelim} we collect some preliminary material
which constitutes the geometric backbone of the paper. The novel
part of the paper is contained in Section \ref{S:proof}. An
essential ingredient in the proof of Theorem \ref{T:permin} is the
second variation formula, see Theorem \ref{T:secondv} below.
Combining the latter with some basic sub-Riemannian integration by
parts formulas, see Lemmas \ref{L:Zf} and \ref{L:Y}, we reduce the
study of the stability of \eqref{ce2i} to checking the validity of
some Hardy type inequalities on $\mS$, see Lemmas \ref{L:X1def} and
\ref{L:nudef}. Using the representation of $\mS$ as a global graph,
we then pull back such Hardy inequalities to ones onto the
$(y,t)-$plane, see Lemma \ref{L:hardyequiv}. Finally, in Lemma
\ref{L:instab} and Corollary \ref{C:instab} we explicitly construct
the directions such that, deforming the surface \eqref{ce2i} along
them, the $H$-perimeter strictly decreases. This establishes the
instability of \eqref{ce2i}.

We would like to close this introduction with some conjectures which
are suggested by the present work. Consider a $C^2$, $H$-minimal
intrinsic $X_1$-graph $\mS\subset \HH$. Denoting by $\mathcal
B_\phi$ the linearized Burger's operator, whose action on a function
$F = F(u,v)$ is given by $\mathcal B_\phi(F)\ =\ F_u + \phi F_v$,
then it was proved in Theorem 1.2 in \cite{ASV} that, provided that
$\phi\in C^2_0(\R^2)$, the $H$-perimeter of $\mS$ is given by
\begin{equation}\label{Pigg}
P_H(\mS)\ =\ \int_{supp(\phi)} \sqrt{1 + \mathcal B_\phi(\phi)^2}\
du\wedge dv\ .
\end{equation}

Now, if we think of \eqref{Pigg} as a functional $P_H(\phi)$ acting
on $\phi$, one can easily recognize that, given $\zeta\in
C^\infty_0(\R^2)$, then the first variation of $P_H(\mS)$ with
respect to the deformation $\mS \to \mS^\la = \mS + \la \mathcal X$,
with $\mathcal X = \zeta X_1$, is given by
\begin{equation}\label{fvig}
 \fv\ \overset{def}{=}\ \frac{d P_H(\mS^\la)}{d\la}\Bl\ =\
\int_{\R^2} \frac{\mathcal B_\phi(\phi)}{{\sqrt{1 + \mathcal
B_\phi(\phi)^2}}} \ (\zeta_u + \phi \zeta_v + \zeta \phi_v)\
du\wedge dv\ .
\end{equation}

In view of Theorem \ref{T:permin}, it is natural to make the
following conjecture: \emph{Suppose that $\phi$, belonging to a
suitable Sobolev space, is a local (or even a global) minimizer of
\eqref{Pigg}, and therefore in particular also a critical point, then
after modification on a set of measure zero, $\phi$ must be of the
type \eqref{vp}}.

If we assume a priori that $\phi\in C^2(\R^2)$, then integrating by
parts in \eqref{fvig} we obtain that
\begin{equation}\label{fvig2}
\fv\ =\ -\ \int_{\R^2} \zeta\ \mathcal B_\phi\left(\frac{\mathcal
B_\phi(\phi)}{\sqrt{1 + \mathcal B_\phi(\phi)^2}}\right)\ du \wedge
dv\ ,\quad\quad \quad\quad \zeta \in C^\infty_0(\R^2)\ .
\end{equation}

On the other hand, under the same regularity hypothesis on $\phi$
one can recognize, see \cite{GS}, that
\begin{equation}\label{MCig}
 \mathcal B_\phi\left(\frac{\mathcal
B_\phi(\phi)}{\sqrt{1 + \mathcal B_\phi(\phi)^2}}\right)\ =\ -\
\mathcal H\ , \end{equation} where $\mathcal H$ represents the
$H$-mean curvature of $\mS$, defined in \eqref{HMCgen} below.
Therefore, a $C^2$ global intrinsic graph $\mS$ is a critical point
of \eqref{Pigg} if and only if $\mS$ is $H$-minimal. As a consequence
of these considerations, in the $C^2$ framework the above conjecture
could be reformulated by saying that: \emph{The only $C^2$, stable,
global intrinsic graphs in $\HH$ are the vertical planes}. Finally,
we would also like to return to the conjecture in the opening of
this introduction and amend it as follows: \emph{In $\HH$ the
vertical planes \eqref{vpi} are the only $C^2$, stable, entire
$H$-minimal graphs (over some plane)}.

\medskip

\noindent \textbf{Acknowledgment:} The problems treated in this
paper were inspired by some stimulating discussions with F. Serra
Cassano and R. Serapioni during a visit of the second named author
at the University of Trento in April 2005. He would like to thank
them for their gracious hospitality. The authors would also like to
thank the anonymous referee for his/her careful reading of the
manuscript and for some comments which helped to improve the
presentation of the paper.

\vskip 0.6in

\section{\textbf{Preliminary material}}\label{S:prelim}

In this section we introduce some relevant notation and
definitions from \cite{DGN3} which will be used in the proof of
Theorem \ref{T:permin}. We consider the first Heisenberg group
$\HH = (\R^3 , \circ)$ with group law
\begin{equation}\label{gl}
 (x,y,t) \ \circ
\ (x',y',t')\ =\ \left(x + x' , y + y' , t + t' + \frac{1}{2} (x
y' - x' y)\right)\ ,
\end{equation}
and non-isotropic dilations $\delta_\la(x,y,t) = (\la x, \la y,
\la^2 t)$, see \cite{S}. Hereafter, we will use the letters $g =
(x,y,t)$, $g' = (x',y',t')$, etc., to indicate points in $\HH$.
Denoting with $(L_g)_*$ the differential of the left-translation
operator $L_g :\HH \to \HH$ defined by $L_g(g') = g \circ g'$, and
letting $e_i$ , $i = 1,2,3$, indicate the standard orthonormal
basis of $\R^3$, one readily verifies that
\begin{equation}\label{p1}
\begin{cases}
X_1(g)\ \overset{def}{=}\ (L_g)_*(e_1)\ =\ \frac{\p}{\p x}\ -\
\frac{y}{2}\ \frac{\p}{\p t}\ ,
\\
X_2(g)\ \overset{def}{=}\ (L_g)_*(e_2)\ =\ \frac{\p}{\p y}\ +\
\frac{x}{2}\ \frac{\p}{\p t}\ ,
\\
T(g)\ \overset{def}{=}\ (L_g)_*(e_3)\ =\ \frac{\p}{\p t}\ .
\end{cases}
\end{equation}

The three vector fields $\{X_1,X_2,T\}$ generate the Lie algebra
$\mathfrak h^1$ of all left-invariant vector fields on $\HH$. They
satisfy at every point of $\HH$ the non-trivial commutation
relation
\begin{equation}\label{p2}
[X_1,X_2]\ =\ T\ , \end{equation} all other commutators being
trivial. In view of \eqref{p2}, the Heisenberg group constitutes the
first (and perhaps most important) prototype of a class of graded
nilpotent Lie groups nowadays known as Carnot groups, see \cite{Fo},
\cite{S}, \cite{Gro}, \cite{Pa2}, \cite{Be}. We observe explicitly
that, if we let $V_1 = \R^2_{x,y} \times \{0\}_t$, and $V_2 =
\{0\}_{x,y} \times \R_t$, then the Heisenberg algebra admits the
decomposition $\mathfrak h^1 = V_1 \oplus V_2$. We assume hereafter
that $\HH$ be endowed with a left-invariant Riemannian metric with
respect to which $\{X_1,X_2,T\}$ constitute an orthonormal basis.
The inner product with respect to this metric will be denoted by
$<\cdot,\cdot>$. This is the only inner product that will be used in
this paper, therefore there will be no confusion with other inner
products, such as for instance the Euclidean one, in $\R^3$. The
corresponding Levi-Civita connection on $\HH$ will be denoted by
$\nabla_X Y$. We will denote by $H\HH$ the subbundle of the tangent
bundle $T\HH$ generated by the distribution $\{X_1,X_2\}$. The
horizontal Levi-Civita connection is given as follows. For any $X\in
\Gamma(T\HH)$, $Y\in \Gamma(H\HH)$ we let
\begin{equation}\label{hlc}
\nabla^H_X Y\ =\ \sum_{i=1}^{2} <\nabla_X Y,X_i> X_i\ ,
\end{equation}
and one can easily verify that $\nabla^H_X Y$ is metric preserving
and torsion free, in the sense that if we define the horizontal
torsion of $\mS$ as
 \[
 T^H(X,Y)\ =\ \nabla^H_X Y\ -\ \nabla^H_Y X\ -\ [X,Y]^H\ ,
 \]
where $[X,Y]^H = \sum_{i=1}^{2} <[X,Y],X_i>X_i$, then $T^H(X,Y)= 0$.
Given a function $f\in C^1(\HH)$, its Riemannian gradient is given
by
\[
\nabla f\ =\ X_1f\ X_1 + X_2f\ X_2 + Tf\ T\ , \] whereas its
horizontal gradient is given by the projection of $\nabla f$ onto
the subbundle $\HH$, hence
\[
\nabh f\ =\ <\nabla f,X_1> X_1 + <\nabla f, X_2> X_2\ =\ X_1f\ X_1 +
X_2f\ X_2\ .
\]

Given an oriented $C^2$ surface $S\subset \HH$, we denote by
$\boldsymbol N$ its (non-unit) Riemannian normal with respect to
the orthonormal frame $\{X_1,X_2,T\}$, and by $\n = \boldsymbol
N/|\boldsymbol N|$ its Riemannian Gauss map. We consider the
quantities
\begin{equation}\label{pq}
p\ =\ <\bN, X_1>\ ,\ \quad q\ =\ <\bN , X_2>\ ,\quad\quad \omega\
=\ <\bN,T>\ ,\ \quad W\ =\ \sqrt{p^2 + q^2}\ .
\end{equation}

In this paper we adopt the classical non-parametric point of view,
see for instance \cite{MM}, according to which a $C^2$ surface
$\mS \subset \HH$ is a subset of $\R^3$ which locally coincides
with the zero set of a real function. Thus, for every $g_0\in \mS$
there exists an open set $\mathcal O\subset \HH$ and a function
$\phi\in C^2(\mathcal O)$ such that: (i) $|\nabla \phi(g)| \not=
0$ for every $g\in \mathcal O$; (ii) $\mS\cap \mathcal O = \{g\in
\mathcal O\mid \phi(g) = 0\}$. We will always assume that $\mS$ is
oriented in such a way that for every $g\in \mS$ one has
\[
\bN(g)\ =\ \nabla \phi(g)\ =\ X_1 \phi(g) X_1\ +\ X_2 \phi(g) X_2\
+\ T \phi(g) T\ . \]

We note explicitly that, in this situation, the functions $p, q,
\omega$ defined in \eqref{pq}, which are given by $p = X_1\phi$, $q
= X_2 \phi$, $\omega = T\phi$, are not only defined on $\mS$, but
for every $g_0\in \mS$ they belong to $C^1(\mathcal O)$. This notion
of $C^2$ surface obviously includes the entire intrinsic graphs
considered in Theorem \ref{T:permin}. In fact, in the case of the
surfaces $\mS$ in \eqref{min}, we have (see \eqref{pqt}),
\[
p = X_1\phi = 1 + \frac{\alpha}{2} y^2\ ,\quad \quad q = X_2 \phi
= - \alpha t - \beta - \frac{\alpha}{2} xy\ ,\quad\quad \omega\ =
T\phi = - \alpha y\ ,
\]
and thus in particular the field $\bN = p X_1 + q X_2 + \omega T$
belongs to $C^\infty(\mathcal O, \HH)$, with $\mathcal O = \HH$.

We emphasize here that the local defining function $\phi$ in (i) and
(ii) above has a different meaning from the function $\phi$ in the
definition of intrinsic graph in the introduction. Given a surface
$\mS \subset \HH$, on the set $\mathcal S\setminus \Sigma(\mS)$ we
define the \emph{horizontal Gauss map} by
\begin{equation}\label{gm}
\nuX\ =\  \pb\ X_1\ +\ \qb\ X_2\ ,
\end{equation}
where we have let
\begin{equation}\label{pbar}
\pb\ =\ \frac{p}{W}\ ,\quad\quad\quad \qb\ =\ \frac{q}{W}\ ,
\quad\quad\text{so that}\quad\quad |\nuX|^2\ =\ \pb^2\ +\ \qb^2\
\equiv\ 1\quad\quad \text{on}\quad\quad \mathcal S \setminus
\Sigma(\mS)\ .
\end{equation}

Given a point $g_0\in \mS\setminus \Sigma(\mS)$, the horizontal
tangent space of $\mS$ at $g_0$ is defined by \[ HT_{g_0}(\mS)\ =\
\{\boldsymbol v \in H_{g_0}\HH\mid <\boldsymbol v,\nuX(g_0)>\ =\
0\}\ .
\]

Let us notice that a basis for $HT_{g_0}(\mS)$ is given by the
vector field
\begin{equation}\label{nup} \nup\ =\ \qb\ X_1\ -\
\pb\ X_2\ .
\end{equation}

Given a function $u\in C^1(\mS)$ one clearly has $\del u(g_0)\in
HT_{g_0}(\mS)$. We next recall some basic definitions from
\cite{DGN3}.

\medskip

\begin{dfn}\label{D:horconS}
Let $\mS\subset \HH$ be a $C^2$ surface, with $\Sigma(\mS) =
\varnothing$, then we define the \emph{horizontal connection} on
$\mS$ as follows. For every $X,Y\in C^1(\mS;HT\mS)$ we let
\[
\delXY\ =\ \nabla^H_{\oX} \oY\ -\ <\nabla^H_{\oX} \oY,\nuX> \nuX\ ,
\]
where $\oX, \oY \in C^1(\HH;H\HH)$ are such that $\oX = X$, $\oY =
Y$ on $\mS$.
\end{dfn}

\medskip

Similarly to the Riemannian case, it is possible to prove that
$\delXY$ does not depend on the extensions $\oX, \oY$. The
tangential horizontal gradient of a function $f\in C^1(\mS)$ is
defined as follows
 \begin{equation}\label{del}
\del f\ =\ \nabh \overline f\ -\ <\nabh \overline f,\nuX> \nuX\ ,
\end{equation}
where $\overline f$ denotes any extension of $f$ to all of $\HH$.
The definition of $\del f$ is well-posed since $\del f$ only depends
on the values of $f$ on $\mS$. Since $|\nuX| \equiv 1$ on $\mS
\setminus \Sigma(\mS)$, we clearly have $<\del f,\nuX> = 0$, and
therefore
\begin{equation}\label{pitagora}
|\del f|^2\ =\ |\nabh \overline f|^2\ -\ <\nabh \overline f,\nuX>^2\
.
\end{equation}

\medskip

All the above definitions are specializations to $\HH$ of analogous
ones for general Carnot groups, see \cite{DGN3}. The next definition
contains the essential geometric concept of horizontal second
fundamental form. It is convenient to state it for $\Hn$, rather
than $\HH$.

\medskip

\begin{dfn}\label{D:sff}
Let $\mS\subset \Hn$ be a $C^2$ hypersurface with $\Sigma(\mS) =
\varnothing$, then for every $X,Y\in C^1(\mS;HT\mS)$ we define a
tensor field of type $(0,2)$ on $\mS$, as follows
\begin{equation}\label{sff} II^{H,\mS}(X,Y)\ =\ <\nabla^H_{\oX} \oY,\nuX>
\nuX\ ,
\end{equation}
where $\oX , \oY$ have the same meaning as in Definition
\ref{D:horconS}. We call $II^{H,\mS}(\cdot,\cdot)$ the
\emph{horizontal second fundamental form} of $\mS$. We also define
$\mathcal A^{H,\mS} : HT \mS \to HT \mS$ by letting for every $g\in
\mS$ and $\boldsymbol u, \boldsymbol v \in HT_{g}\mS$
\begin{equation}\label{shape}
<\mathcal A^{H,\mS} \boldsymbol u,\boldsymbol v>\ =\ -\
<II^{H,\mS}(\boldsymbol u,\boldsymbol v),\nuX>\ =\ -\
<\nabla_{\oX}^H \oY,\nuX>\ ,
\end{equation}
where $X, Y \in C^1(\mS,HT \mS)$ are such that $X_g = \boldsymbol
u$, $Y_g = \boldsymbol v$, and $\oX, \oY$ are as above. We call the
linear map $\mathcal A^{H,\mS} : HT_{g}\mS \to HT_{g}\mS$ the
\emph{horizontal shape operator}. If $\boldsymbol
e_1,...,\boldsymbol e_{2n-1}$ denotes a local orthonormal frame for
$HT\mS$, then the matrix of the horizontal shape operator with
respect to the basis $\boldsymbol e_1,...,\boldsymbol e_{2n-1}$ is
given by the $(2n-1)\times(2n-1)$ matrix \ $-
\big[<\nabla_{\boldsymbol e_i}^H \boldsymbol
e_j,\nuX>\big]_{i,j=1,...,2n-1}$.

\end{dfn}

\medskip

If $\mS$ has non-empty characteristic locus $\Sigma(\mS)$, then we
consider $\mS' = \mS \setminus \Sigma(\mS)$ and define the
$H$-\emph{mean curvature} of $\mS$ at a point $g_0\in \mS'$ as
follows
\begin{equation}\label{HMCgen}
\mathcal H\ =\ -\ trace\ \mathcal A^{H,\mS}\ =\ -\ \sum_{j=1}^{2n-1}
<\nabla_{\boldsymbol e_i}^H \boldsymbol e_j,\nuX>\ .
\end{equation}

We recall that is was proved in \cite{B}, \cite{Ma3} that $\mathcal
H^{Q-1}(\Sigma(\mS)) = 0$, where $\mathcal H^s$ denotes the
$s$-dimensional Hausdorff measure associated with the horizontal or
Carnot-Carath\'eodory distance of $\bG$, and $Q$ indicates the
homogeneous dimension of $\bG$. If $g_0\in \Sigma(\mS)$ we let
\[
\mathcal H(g_0)\ =\ \underset{g\to g_0, g\in \mathcal S\setminus
\Sigma(\mS)}{\lim}\ \mathcal H(g)\ ,
\]
provided that such limit exists, finite or infinite. We do not
define the $H$-mean curvature at those points $g_0\in \Sigma(\mS)$
at which the limit does not exist. The following result is taken
from \cite{DGN3}.

\medskip

\begin{prop}\label{P:equalMC}
The $H$-mean curvature of $\mS \subset \HH$ coincides with the
function
\begin{equation}\label{equal2}
\mathcal H\ =\ \sum_{i=1}^{2}\ \di\ <\nuX,X_i>\ =\ \duno \pb + \ddue
\qb\ =\ X_1 \pb\ +\ X_2 \qb\ ,
\end{equation}
where $\pb , \qb$ are as in \eqref{pbar}.
\end{prop}

\medskip

\begin{dfn}\label{D:minimal}
A $C^2$ surface $\mS\subset \HH$ is called $H$-minimal if
$\mathcal H \equiv 0$ as a continuous function on $\mS$.
\end{dfn}

\medskip

In keeping up with the notation of \cite{DGN3} it will be convenient
to indicate with $Y\zeta$ and $Z\zeta$ the respective actions of the
vector fields $\nuX$ and $\nup$ on a function $\zeta\in
C^1_0(\mathcal S \setminus \Sigma(\mS))$, thus
\begin{equation}\label{Y}
Y\zeta\ \overset{def}{=}\ \pb\ X_1 \zeta\ +\ \qb\ X_2 \zeta\
,\quad\quad\quad Z\zeta\ \overset{def}{=}\ \qb\ X_1\zeta\ -\ \pb\
X_2\zeta\ .
\end{equation}

The frame $\{Z,Y,T\}$ is orthonormal. It is worth observing that,
since the metric tensor $\{g_{ij}\}$ with respect to the inner
product $<\cdot,\cdot>$ has the property $det\{g_{ij}\} = 1$, then
the (Riemannian) divergence in $\HH$ of these vector fields is
given by
\begin{equation}\label{divY}
div\ Y\ =\ X_1 \pb\ +\ X_2 \qb\ =\ \mathcal H\ ,\quad\quad\quad
div\ Z\ =\ X_1 \qb\ -\ X_2 \pb\ .
\end{equation}

Using Cramer's rule one easily obtains from \eqref{Y}
\begin{equation}\label{Xs}
X_1 \zeta\ =\ \qb\ Z\zeta\ +\ \pb\ Y\zeta\ ,\quad\quad\quad
X_2\zeta\ =\ \qb\ Y\zeta\ -\ \pb\ Z\zeta\ .
\end{equation}

One also has
\begin{equation}\label{deltas}
\duno \zeta\ =\ \qb\ Z\zeta\ ,\quad\quad\quad\quad \ddue \zeta\ =\
-\ \pb\ Z\zeta\ ,
\end{equation}
so that
\begin{equation}\label{delZ}
|\del \zeta|^2\ =\ (Z\zeta)^2\ .
\end{equation}

We notice that
\begin{equation}\label{mc}
\qb Z\pb\ -\ \pb Z\qb\ =\ \mathcal H\ . \end{equation}

This can be easily recognized using Proposition \ref{P:equalMC} and
\eqref{Xs}, as follows
\[ \mathcal H\ =\ X_1
\pb\ +\ X_2 \qb\ =\ \qb Z\pb - \pb Z\qb + \pb Y\pb + \qb Y\qb\ =\
\qb Z\pb - \pb Z\qb\ ,
\]
where we have used the fact that $0 = \frac{1}{2} Y(\pb^2 + \qb^2) =
\pb Y\pb + \qb Y\qb$. Finally, we will need the following identity
\begin{equation}\label{accasquare}
\mathcal H^2\ =\ (Z\pb)^2\ +\ (Z\qb)^2\ .
\end{equation}

This can be easily proved observing that \eqref{mc} and the identity
$\pb^2 + \qb^2 = 1$ give
\[
\mathcal H^2\ =\ (Z\pb)^2 + (Z\qb^2) - (\pb Z\pb + \qb Z\qb)\ =\
(Z\pb)^2 + (Z\qb^2)\ .
\]

In the classical theory of minimal surfaces, the concept of area or
perimeter occupies a central position, see \cite{DG1}, \cite{DG2},
\cite{DCP}, \cite{G}, \cite{MM}. In sub-Riemannian geometry there
exists an appropriate notion of perimeter. Given an open set $\Om
\subset \HH$ we denote $\mathcal F(\Om) = \{\zeta \in
C^1_0(\Om;H\HH)\mid ||\zeta||_{L^\infty(\Om)}\leq 1\}$. A function
$u\in L^1(\Om)$ is said to belong to $BV_H(\Om)$ (the space of
functions with finite horizontal bounded variation), if
\[
Var_H(u;\Om)\ =\ \underset{\zeta\in \mathcal F(\Om)}{\sup}\
\int_\Om u\ div_H \zeta\ dg\ < \ \infty\ .
\]

This space becomes a Banach space with the norm $||u||_{BV_H(\Om)} =
||u||_{L^1(\Om)} + Var_H(u;\Om)$. Given a measurable set $E\subset
\HH$, the $H$-\emph{perimeter} of $E$ with respect to the open set
$\Om\subset \HH$ is defined as follows, see for instance \cite{CDG},
and \cite{GN},
\[
P_H(E;\Om)\ =\ Var_H(\chi_E;\Om)\ .
\]

Given an oriented $C^2$ surface $\mS$, we will denote by $d\sigma_H$
the $H$-perimeter measure concentrated on $\mS$. For any Borel
subset $E\subset \mS$ such that $P_H(E)<\infty$, one has
\begin{align}\label{sigmah}
P_H(E)\ & =\ \int_E d\sigma_H\ =\ \int_E \sqrt{<\n,X_1>^2 + <\n ,
X_2>^2}\ d\sigma
\\
& =\ \int_E \frac{\sqrt{<\bN,X_1>^2 + <\bN , X_2>^2}}{|\bN|}\
d\sigma\ =\ \int_E \frac{W}{|\bN|}\ d\sigma\ , \notag
\end{align}
where in the last equality we have used \eqref{pq}. We thus obtain
from \eqref{sigmah}
\begin{equation}\label{pm}
d\sigma_H\ =\ \frac{W}{|\bN|}\ d\sigma\ ,
\end{equation}
where $d\sigma$ denotes the standard surface measure.

\vskip 0.6in


\section{\textbf{Proof of Theorem \ref{T:permin}}}\label{S:proof}

\vskip 0.2in

This section is devoted to proving Theorem \ref{T:permin}. In the
course of the proof we need to build some auxiliary results which we
fell have an independent interest. We begin by recalling the notions
of first and second variation of the $H$-perimeter introduced in
Definition \ref{D:minimali}. Classical minimal surfaces are critical
points of the perimeter (area functional). It is natural to ask what
is the connection between the notion of $H$-minimal surface and that
of $H$-perimeter. The answer to this question is contained in the
following result from \cite{DGN3}, see also \cite{DGN2}.

\medskip

 \begin{thrm}\label{T:variations}
Let $\mS\subset \HH$ be an oriented $C^2$ surface with $\Sigma(\mS)
= \varnothing$, then
\begin{equation}\label{fvH}
\fv\ \ =\
 \int_{\mathcal S}
\mathcal H\ \frac{\cos(\mathcal X \angle \bN)}{\cos(\nuX \angle
\bN)}\ |\mathcal X|\ d\sigma_H\ ,
\end{equation}
where $\angle$ denotes the angle between vectors in the inner
product $<\cdot,\cdot>$. In particular, $\mathcal S$ is stationary
if and only if it is $H$-minimal.
\end{thrm}

\medskip

We emphasize that, thanks to the assumption $\Sigma(\mS) =
\varnothing$, the denominator in the integrand in the right-hand
side of \eqref{fvH} does not vanish on $\mS$. We mention that
versions of Theorem \ref{T:variations} have also been obtained
independently by other people. An approach based on motion by
$H$-mean curvature can be found in \cite{BC}. When $\mathcal X = h
\nuX$, then a proof based on CR-geometry can be found in
\cite{CHMY}.

\medskip

A central (and more complex) result for this paper is the following
theorem established in \cite{DGN3}. Recalling the function $\omega$
defined in \eqref{pq}, henceforth we let $\ob = \omega/W$.

\medskip

\begin{thrm}\label{T:secondv}
Let $\mS\subset \HH$ be a $C^2$ oriented surface with $\Sigma(\mS) =
\varnothing$. The second variation of the $H$-perimeter with respect
to the deformation $\mS \to \mS^\la = \mS + \la \mathcal X$, with
$\mathcal X = a X_1 + b X_2 + k T\in C^2_0(\mS;\HH)$, is given by
the formula
\begin{align}\label{svH}
\sv\ & =\ \int_{\mathcal S} \bigg\{ -\ 2\ (\pb Zb - \qb Za) \left(Tk
- \ob Yk\right)
\\
& +\ \left(Ta - \ob Ya\right) \bigg[- 2 \qb Zk - \qb
(a \pb + b \qb) - \pb (a \qb - b \pb)\bigg] \notag\\
& +\ \left(Tb - \ob Yb\right)\bigg[2 \pb Zk + \pb (a
\pb + b \qb) - \qb (a \qb - b \pb)\bigg] \notag\\
& +\ 2\ (a \qb - b \pb) (\qb Za - \pb Zb)\ \ob
\notag\\
& +\ \left(Za + \ob\ \pb\ Zk\right)^2\ +\ \left(Zb + \ob\ \qb\
Zk\right)^2
\notag\\
& +\  (a^2 + b^2)\ \ob^2
\notag\\
& +\ 2\ \ob (a Za + b Zb)\ +\ 2\
\ob^2 (a \pb + b \qb) Zk \notag\\
& -\ \left(\qb Za - \pb Zb + (a \qb  - b \pb) \ob\right)^2\bigg\}\
d\sigma_H . \notag
\end{align}
\end{thrm}

\medskip

\begin{cor}\label{C:X1def}
If we choose $a\in C^\infty_0(\mS)$, $b \equiv k\equiv 0$, and
therefore $\mathcal X = a X_1$, then the corresponding second
variation of the $H$-perimeter is given by
\begin{align}\label{X1def}
\sv\ & =\ \int_{\mathcal S} \bigg\{\pb^2 (Za)^2 + \pb^2\ \ob^2\ a^2
\\
& +\ \ob Z(a^2) - \pb\ \qb \left(T(a^2) - \ob
Y(a^2)\right)\bigg\}\ d\sigma_H\ . \notag
\end{align}
\end{cor}

\begin{proof}[\textbf{Proof}]
It follows in an elementary fashion from \eqref{svH}. One only
needs to keep in mind that $\pb^2 + \qb^2 = 1$.

\end{proof}

\medskip

\begin{cor}\label{C:nuXdef}
Given a $C^3$ oriented surface $\mS \subset \HH$, with $\Sigma(\mS)
= \varnothing$, consider the deformation
\begin{equation}\label{nuXdef} \mS^\la\ =\ \mS\ +\ \la\ (h\
\nuX\ +\ k\ T)\ , \quad\quad\quad h , k \in C^2_0(\mS)\
,\end{equation} corresponding to the choice $\mathcal X = \pb h X_1
+ \qb h X_2$ (notice that $a= \pb h , b = \qb h \in C^2_0(\mS)$).
One has
\begin{align}\label{nuXdef}
\sv\ & =\ \int_{\mathcal S} \bigg(Zh + \ob Zk\bigg)^2\ d\sigma_H
\\
& +\ 2\ \int_\mS h\ \mathcal H\ \left(Tk - \ob
Yk\right)\ d\sigma_H \notag\\
& +\ \int_\mS \bigg\{\ob Z(h^2) + 2\ \mathcal A\ h\ Zk\ +\
\mathcal A\ h^2\bigg\}\ d\sigma_H , \notag
\end{align}
where we have set
\begin{equation}\label{alpha}
\mathcal A\ =\ (\pb T\qb - \qb T\pb)\ +\ \ob(\qb Y\pb - \pb Y\qb)\
+\ \ob^2\ .
\end{equation}
\end{cor}

\begin{proof}[\textbf{Proof}]
We notice that we presently have
\[
\pb a + \qb b\ =\ h\ ,\quad\quad\quad \qb a - \pb b\ =\ 0\ ,
\]
\[
a\ Za\ +\ b\ Zb\ =\ h\ Zh\ ,
\]
\[
\pb Zb - \qb Za\ =\ h (\pb Z\qb - \qb Z\pb)\ =\ -\ h\ \mathcal H\
,
\]
where in the last equality we have used the identity \eqref{mc}.
We also have

\[
Ta - \ob Ya\ =\ \pb \left(Th - \ob Yh\right)\ +\ h\ \left(T\pb -
\ob Y\pb\right)\ ,
\]
\[
Tb - \ob Yb\ =\ \qb \left(Th - \ob Yh\right)\ +\ h\ \left(T\qb -
\ob Y\qb\right)\ ,
\]
\begin{align*}
\left(Za + \ob\ \pb\ Zk\right)^2\ & =\ \pb^2 \left(Zh + \ob
Zk\right)^2
\\
& +\ (Z\pb)^2\ h^2\ +\ 2\ \pb\ Z\pb\ h\ \left(Zh + \ob Zk\right)\
,
\end{align*}
\begin{align*}
\left(Zb + \ob\ \qb\ Zk\right)^2\ & =\ \qb^2 \left(Zh + \ob
Zk\right)^2
\\
& +\ (Z\qb)^2\ h^2\ +\ 2\ \qb\ Z\qb\ h\ \left(Zh + \ob Zk\right)\
.
\end{align*}

We next observe that

Substituting these formulas in the right-hand side of \eqref{svH},
and using \eqref{accasquare}, we reach the desired conclusion.

\end{proof}

\medskip

A different approach to a version of \eqref{nuXdef} based on
CR-geometry was found in \cite{CHMY}. To reduce further the
expressions in the right-hand side of \eqref{X1def},
\eqref{nuXdef} we would like to transform the terms containing the
derivatives $Z(a^2)$, $T(a^2)$, $Y(a^2)$, $Z(h^2)$. For this, we
will use the following basic integration by parts formulas proved
in \cite{DGN3}.

\medskip

\begin{lemma}\label{L:Zf}
Let $\mS\subset \HH$ be a $C^2$ oriented surface with $\Sigma(\mS) =
\varnothing$. For any $\zeta\in C_0^1(\mS)$ one has
\[
\int_{\mathcal S} Z\zeta\ d\sigma_H\ =\ -\ \int_{\mathcal S}
\zeta\ \ob\  d\sigma_H\ .
\]
\end{lemma}

\medskip

\begin{lemma}\label{L:Y}
With $\mS$ as in Lemma \ref{L:Zf}, for any function $\zeta\in
C^1_0(\mathcal S)$ one has
\[
 \int_\mathcal S T\zeta\ d\sigma_H\ =\ \int_\mathcal S Y\zeta\
 \ob\ d\sigma_H\ +\ \int_\mathcal S \zeta\
 \ob\ \mathcal H\  d\sigma_H\ .
 \]
\end{lemma}

\medskip

Using Lemma \ref{L:Zf}  we find
\begin{align}\label{Zasquare}
& \int_\mS \ob Z(a^2)\ d\sigma_H\ =\  \int_\mS Z \left(\ob
a^2\right) d\sigma_H\ -\ \int_\mS a^2 Z\ob\ d\sigma_H
\\
& =\ -\ \int_\mS a^2 \ob^2\ d\sigma_H\ -\ \int_\mS a^2 Z\ob\
d\sigma_H\ . \notag
\end{align}

From Lemma \ref{L:Y} we obtain instead
\begin{align}\label{Tasquare}
& -\ \int_\mS  \pb\ \qb \left(T(a^2) - \ob
Y(a^2)\right)\ d\sigma_H\\
& =\ \int_\mS a^2 \bigg\{(\pb T\qb + \qb T\pb) - \ob(\pb Y\qb +
\qb Y\pb) - \pb\ \qb\ \ob \mathcal H\bigg\}\ d\sigma_H\ . \notag
\end{align}

Substituting \eqref{Zasquare}, \eqref{Tasquare} into
\eqref{X1def}, and keeping \eqref{del} in mind, we finally obtain.

\medskip

\begin{lemma}\label{L:X1def}
Let $\mS\subset \HH$ be a $C^2$, oriented surface, with $\Sigma(\mS)
= \varnothing$, then the second variation of the $H$-perimeter, with
respect to deformation $\mS^\la = \mS + \la a X_1$, is given by
\begin{align}\label{2varX1}
\mathcal V^H_{II}(\mS;aX_1)\ & =\ \int_{\mathcal S} \pb^2 |\del a|^2
d\sigma_H\ +\ \int_\mS a^2 \bigg\{(\pb T\qb + \qb T\pb) - \ob (\pb
Y\qb + \qb Y\pb)
\\
& -\ \qb^2 \ob^2 - Z \ob - \pb\ \qb\ \ob \mathcal H\bigg\}\
d\sigma_H\ . \notag
\end{align}
\end{lemma}

\medskip

To establish the next lemma we need the following auxiliary result.

\medskip

\begin{lemma}\label{L:Zob}
On a surface $\mS$ as in Lemma \ref{L:X1def}, one has \[
-\ Z\ob\ =\
\mathcal A\ ,
\]
where $\mathcal A$ is the quantity defined in \eqref{alpha}.
\end{lemma}

\begin{proof}[\textbf{Proof}]
From the definition of $\ob$ one has
\begin{equation}\label{Zob1}
-\ Z \ob\ =\ \ob\ \frac{ZW}{W}\ -\ \frac{Z\omega}{W}\ .
\end{equation}

We now claim that
\begin{equation}\label{Zob2}
\frac{Z\omega}{W}\ =\ \qb\ T \pb\ -\ \pb\ T \qb\ ,
\end{equation}
and that, furthermore,
\begin{equation}\label{Zob3}
\frac{ZW}{W}\ =\ \qb\ Y\pb\ -\ \pb\ Y\qb\ +\ \ob\ .
\end{equation}

It should be obvious to the reader that, inserting \eqref{Zob2},
\eqref{Zob3} into \eqref{Zob1}, we obtain the desired conclusion. We
are thus left with proving \eqref{Zob2} and \eqref{Zob3}. For the
former, we observe that
\[
<Z,\bN>\ =\ 0\ .
\]

If $\phi$ denotes a local defining function of $\mS$ in the
neighborhood of an arbitrary point, we thus have $Z\phi = 0$.
Applying $T$ to this identity, we obtain
\begin{align*}
0 \ =\ T(Z\phi) &\ =\ T(\qb X_1\phi\ -\ \pb X_2\phi)
\ =\ T\qb X_1\phi\ +\ \qb\, TX_1\phi\ -\ T\pb X_2\phi\ -\ \pb\, TX_2\phi\\
& \ =\ T\qb\,X_1\phi - T\pb\,X_2\phi\ +\ \qb X_1T\phi - \pb X_2T\phi
\ =\ pT\qb - q T\pb + Z(T\phi)\ ,
\end{align*}
where we have used $[X_i,T]=0$, $i =1,2$. It follows that
\[
\frac{Z\omega}{W}\ =\ \frac{Z(T\phi)}{W}\ =\ \qb\, T\pb - \pb\,
T\qb\ ,
\]
which proves \eqref{Zob2}. As for \eqref{Zob3}, we  have
\begin{equation*}
\omega\ =\ T\phi\ =\ X_1X_2\phi\ -\ X_2X_1\phi\ =\ X_1(\qb\ W)\ -\
X_2(\pb\ W)\ =\  -\ (X_2\pb\ -\ X_1\qb)\ W\ +\ ZW\
 , \end{equation*}
from which the desired conclusion follows immediately.

\end{proof}

\medskip

Using Lemma \ref{L:Zf} in \eqref{nuXdef} of Corollary
\ref{C:nuXdef}, in combination with Lemma \ref{L:Zob}, we obtain.

\medskip

\begin{lemma}\label{L:nudef}
Let $\mS\subset \HH$ be a $C^3$, oriented surface, with $\Sigma(\mS)
= \varnothing$, then the second variation of the $H$-perimeter with
respect to the deformation of $\mS^\la = \mS + \la h \nuX$, is given
by
\begin{align}\label{2varnu}
\mathcal V_{II}^H(\mS;h\nuX)\ & =\ \int_\mS \bigg\{(Zh)^2 +
h^2\left(2 \mathcal A - \ob^2\right)\bigg\}\ d\sigma_H
\\
& =\ \int_\mS \bigg\{(Zh)^2 + h^2 \big[2 (\pb T\qb - \qb T\pb) + 2
\ob (\qb Y\pb - \pb Y\qb) + \ob^2\big]\bigg\}\ d\sigma_H\ , \notag
\end{align}
where $\mathcal A$ is defined in \eqref{alpha}.
\end{lemma}

\medskip

After these preparations we turn to the core of the proof of
Theorem \ref{T:permin}. We will focus on the case in which the
surface is given by
\begin{equation}\label{ce2}
\mS\ =\ \{(x,y,t)\in \HH \mid x = y (\alpha t + \beta)\}\ ,
\end{equation}
the other family of surfaces $\mS = \{(x,y,t)\in \HH \mid y = x
(\alpha t + \beta)\}$, with $\alpha < 0$ and $\beta\in \R$, being
treated by completely analogous considerations. Our first step in
the proof of Theorem \ref{T:permin} is to compute the second
variation of the $H$-perimeter for $\mS$. In view of
\eqref{2varX1} in Lemma \ref{L:X1def}, or \eqref{2varnu} in Lemma
\ref{L:nudef}, we need to compute the quantities which appear as
the coefficient of $a^2$ and $h^2$ in the integral in the
right-hand side of the respective formulas. This is the content of
the next lemma.

\medskip

\begin{lemma}\label{L:coefficient}
Let $\mS$ be the $H$-minimal surface given by \eqref{ce2}, then
one has
\begin{align}\label{1def}
& (\pb T\qb + \qb T\pb) - \ob (\pb Y\qb + \qb Y\pb)
-\ \qb^2 \ob^2 \\
& - Z \ob - \pb\ \qb \ \ob\ \mathcal H\ =\ -\ \frac{2 \alpha}{W^2
(1 + (\alpha t + \beta)^2)}\ , \notag
\end{align}
\begin{align}\label{2def}
& 2 (\pb T\qb - \qb T\pb) + 2 \ob(\qb Y\pb - \pb Y\qb) + \ob^2 \
=\ -\ \frac{2 \alpha}{W^2}\ .
\end{align}
\end{lemma}

\begin{proof}[\textbf{Proof}]
We can use the global defining function $\phi(x,y,t) = x - y
(\alpha t + \beta)$. As previously stipulated, we assume that
$\mS$ is oriented in such a way that $\bN = \nabla \phi = X_1\phi\
X_1 + X_2\phi\ X_2 + T\phi\ T$. Recalling \eqref{pq}, simple
calculations based on \eqref{p1} thus give
\begin{equation}\label{pqt}
p = X_1\phi = 1 + \frac{\alpha}{2} y^2\ ,\quad \quad q = X_2 \phi
= - \alpha t - \beta - \frac{\alpha}{2} xy\ ,\quad\quad \omega\ =
T\phi = - \alpha y\ .
\end{equation}

The second equation in \eqref{pqt} becomes on $\mS$
\begin{equation}\label{qons}
q\ =\ -\ (\alpha t + \beta) \left(1\ +\
\frac{\alpha}{2}y^2\right)\ .
\end{equation}

We thus find on $\mS$
\begin{equation}\label{Wons}
W^2\ =\ |\nabh \phi|^2\ =\ \left(1\ +\
\frac{\alpha}{2}y^2\right)^2\ (1 + (\alpha t + \beta)^2)\ .
\end{equation}

Using \eqref{Y} we obtain
\begin{equation}\label{ztf}
 Z(T\phi)\ =\
\qb\ X_1(T\phi)\ -\ \pb\ X_2(T\phi)\ =\ \frac{\alpha p}{W}\ =\
\frac{\alpha\left(1 + \frac{\alpha}{2}y^2\right)}{W}\ >\ 0\ .
\end{equation}

Next, we have
\begin{equation}\label{hardy1}
X_1p\ =\ 0\ ,\quad\quad X_1 q\ =\ 0\ ,\quad\quad X_2p\ =\ \alpha
y\ ,\quad\quad X_2q\ =\ -\ \alpha x\ .
\end{equation}

This gives
\begin{equation}\label{XW}
X_1 W\ =\ \frac{p X_1p + q X_1q}{W}\ =\ 0\ ,\quad\quad X_2 W\ =\
\frac{p X_2p + q X_2q}{W}\ =\ \alpha y \left(1 + (\alpha t +
\beta)^2)\right)^{1/2}\ .
\end{equation}

From \eqref{Y} and \eqref{XW} we find
\begin{equation}\label{ZW}
ZW\ =\ \qb X_1W - \pb X_2 W\ =\ -\ \frac{p}{W} X_2W\ =\ -\ \alpha\
y\ .
\end{equation}

Combining \eqref{ZW} with \eqref{ztf} we obtain
\begin{equation}\label{quotient}
Z \ob\ =\ \frac{Z(T\phi)}{W} - \frac{T\phi}{W^2} ZW\ =\
\frac{\alpha - \frac{\alpha^2}{2} y^2}{W^2}\ .
\end{equation}

Using the above formulas it is not difficult to verify that
\begin{equation}\label{Ypq}
Y\pb\ =\ Y\qb\ =\ 0\ .
\end{equation}

We now have from \eqref{pqt} \[ TW\ =\ \frac{p Tp + q Tq}{W}\ =\
-\ \frac{\alpha q}{W}\ ,
\]
and therefore we easily find
\begin{equation}\label{Tpq}
\begin{cases}
T\pb\ =\ -\ \frac{\alpha(\alpha t + \beta)}{\left(1 +
\frac{\alpha}{2}y^2\right) (1 + (\alpha t + \beta)^2)^{3/2}}\ =\
-\ \frac{\alpha(\alpha t + \beta)}{W (1 + (\alpha t + \beta)^2)}\
,
\\
T\qb\ =\ -\ \frac{\alpha}{W (1 + (\alpha t + \beta)^2)}\ .
\end{cases}
\end{equation}

From \eqref{pqt} and \eqref{Tpq} we conclude that
\begin{equation}\label{TpTq}
\pb\ T\qb\ +\ \qb\ T\pb \ =\ \frac{\alpha \left(1 +
\frac{\alpha}{2}y^2\right)((\alpha t + \beta)^2 - 1)}{W^2 (1 +
(\alpha t + \beta)^2)}\ .
\end{equation}

From \eqref{qons}, \eqref{quotient}, and \eqref{TpTq}, and
elementary computations, we easily reach the conclusion that
\eqref{1def} holds. In a similar fashion, we obtain the proof of
\eqref{2def} by \eqref{pqt}, \eqref{qons}, \eqref{Wons} and
\eqref{Tpq}.

\end{proof}

\medskip

From Lemmas \ref{L:X1def}, \ref{L:nudef} and \ref{L:coefficient},
we obtain the following corollary.

\medskip

\begin{cor}\label{C:X1graph}
Let $\mS$ be the $H$-minimal surface given by \eqref{ce2}. For any
$a\in C^\infty_0(\mS)$, the second variation along the deformation
$\mS \to \mS^\la = \mS + \la a X_1$ is given by
\begin{align}\label{2varX1graph}
\mathcal V^H_{II}(\mS;aX_1)\ & =\ \int_{\mathcal S} \frac{\left(1 +
\frac{\alpha}{2}y^2\right)^2}{W^2} |\del a|^2  d\sigma_H\ -\ 2\
\alpha\ \ \int_\mS \frac{a^2}{W^2 (1 + (\alpha t + \beta)^2)}\
d\sigma_H\ .
\end{align}
For any $h\in C^\infty_0(\mS)$, the second variation along the
deformation $\mS \to \mS^\la = \mS + \la h \nuX$ is given by
\begin{align}\label{2varnugraph}
\mathcal V^H_{II}(\mS;h \nuX)\ & =\ \int_{\mathcal S} |\del h|^2
d\sigma_H\ -\ 2\ \alpha\ \ \int_\mS \frac{h^2}{W^2}\ d\sigma_H\ .
\end{align}
\end{cor}

\medskip

We now consider the global smooth parametrization $\theta :
\mathbb R^2 \to \mathbb R^3$ of the surface $\mS$ given by
$\theta(y,t) = (y(\alpha t + \beta),y,t)$.  Clearly, $\mS =
\theta(\mathbb R^2)$.

\medskip

\begin{lemma}\label{L:hardyequiv}
Let $\mS$ be the $H$-minimal surface given by \eqref{ce2}. For any
$a\in C^\infty_0(\mS)$,  then one has
\begin{align}\label{he}
\mathcal V^H_{II}(\mS;aX_1)\ & =\ \int_{\mathbb R^2} \frac{\left(1 +
\frac{\alpha}{2}y^2\right)\ u_y^2}{(1 + (\alpha t +
\beta)^2)^{3/2}}\ dy dt\\
& -\ 2 \alpha \int_{\mathbb R^2} \frac{u^2}{\left(1 +
\frac{\alpha}{2}y^2\right) (1 + (\alpha t + \beta)^2)^{3/2}}\ dy
dt\ ,\notag
\end{align}
where $u = a \circ \theta \in C^\infty_0(\mathbb R^2)$. For any
$h\in C^\infty_0(\mS)$, the one has
\begin{align}\label{he2}
\mathcal V^H_{II}(\mS;h \nuX)\ & =\ \int_{\mathbb R^2} \frac{\left(1
+ \frac{\alpha}{2}y^2\right)\ u_y^2}{(1 + (\alpha t +
\beta)^2)^{1/2}}\ dy dt\\
& -\ 2 \alpha \int_{\mathbb R^2} \frac{u^2}{\left(1 +
\frac{\alpha}{2}y^2\right) (1 + (\alpha t + \beta)^2)^{1/2}}\ dy
dt\ ,\notag
\end{align}
where this time we have set $u = h \circ \theta \in
C^\infty_0(\mathbb R^2)$.
\end{lemma}

\begin{proof}[\textbf{Proof}]
In order to prove \eqref{he} we make some reductions. Keeping in
mind that from \eqref{pm} we have $d\sigma_H = (|\nabh
\phi|/|\nabla \phi|) d\sigma = (W/|\nabla \phi|) d\sigma$, from
\eqref{Wons} we obtain
\[
\int_\mS \frac{a^2}{W^2 (1 + (\alpha t + \beta)^2)} d\sigma_H\  =\
\int_{\mathbb R^2} \frac{u^2}{\left(1\ +\
\frac{\alpha}{2}y^2\right) (1 + (\alpha t + \beta)^2)^{3/2}}\ dy
dt\ .
\]

In order to express the first integral in the right-hand side of
\eqref{2varX1graph} as an integral on $\mathbb R^2$, we compute
$|\del a|^2$. We have from \eqref{del}, \eqref{Y}, \eqref{pqt} and
\eqref{qons}
\begin{align}\label{delh}
|\del a|^2\ & =\ (Za)^2\ =\ (\qb X_1 a - \pb X_2 a)^2
\\
& =\ \frac{((\alpha t + \beta) X_1a + X_2 a)^2}{1 + (\alpha t +
\beta)^2}\ . \notag
\end{align}

Now, the chain rule gives $u_y = (\alpha t + \beta)a_x + a_y$, and
therefore we see from \eqref{p1} that we have on $\mS$
\[
(\alpha t + \beta) X_1a + X_2 a\ =\ (\alpha t + \beta)\ a_x + a_y\
=\ u_y\ .
\]

From \eqref{delh} we thus conclude that
\[
\int_{\mathcal S} \pb^2 \left|\delta_H a\right|^2 d\sigma_H\ =\
\int_{\mathbb R^2} \frac{\left(1\ +\ \frac{\alpha}{2}y^2\right)\
u_y^2}{(1 + (\alpha t + \beta)^2)^{3/2}}\ dy dt\ .
\]

This proves \eqref{he}. The proof of \eqref{he2} proceeds
analogously, and we omit the details.

\end{proof}

\medskip

\begin{lemma}\label{L:instab}
Let $\chi_k \in C^\infty_0(\R)$ be such that $0\leq \chi_k(s) \leq
1$, $\chi_k(s) = 0$ for $|s| > 2k$, $\chi_k(s) \equiv 1$ for $|s|<
k$, and $|\chi_k'(s)| \leq C/k$ with $C$ independent of $k$.
Define for any $\alpha >0$
\[
f_k(y)\ = \ \frac{\chi_k(y)}{\sqrt{1 + \frac{\alpha}{2} y^2}} \ .
\]
We have $f_k\in C^\infty_0(\R)$, and there exists $k_0 \in \mathbb
N$ such that for all $k > k_0$
\begin{equation}\label{reverse_ineq}
\int_{\R} \frac{f_k(y)^2}{1 + \frac{\alpha}{2}y^2}\,dy \ > \
\frac{1}{2\alpha} \int_{\R} \left(1 + \frac{\alpha}{2}y^2\right)
f_k'(y)^2\,dy \ .
\end{equation}
\end{lemma}

\begin{proof}[\textbf{Proof}]
We begin by observing that integration of the right-hand side of
\eqref{reverse_ineq} gives
\begin{align}\label{rhs}
& \int_{\R} (1 + \frac{\alpha}{2}y^2) f_k'(y)^2\,dy \\
\notag &\quad\ =\  \int_\R (\chi_k')^2\ dy\ +\ \frac{\alpha^2}{4}
\int_\R \frac{y^2 \chi_k^2}{(1 + \frac{\alpha}{2}y^2)^2} dy\ +\
\frac{\alpha}{2} \int_\R \chi_k^2 \left(\frac{y}{1 +
\frac{\alpha}{2}y^2}\right)' dy\ . \notag
\end{align}

Observing that
\[
\left(\frac{y}{1 + \frac{\alpha}{2}y^2}\right)' \ =\ \frac{1 -
\frac{\alpha}{2} y^2}{(1 + \frac{\alpha}{2}y^2)^2}\ ,
\]
we conclude from \eqref{rhs} and dominated convergence
\begin{align}\label{rhs2}
& \int_{\R} (1 + \frac{\alpha}{2}y^2) f_k'(y)^2\,dy\ =\
O\left(\frac{1}{k}\right)\  +\ \frac{\alpha}{2} \int_{\R}
\frac{\chi_k^2}{\left(1 + \frac{\alpha}{2}y^2\right)^2}\ dy
\\
&  \longrightarrow\ \frac{\alpha}{2} \int_{\R} \frac{1}{\left(1 +
\frac{\alpha}{2}y^2\right)^2}\ dy\ =\ \sqrt{\frac{\alpha}{2}}\
\frac{\pi}{2}\ ,\quad\quad \text{as}\quad k \to \infty\ . \notag
\end{align}

On the other hand, again by dominated convergence, we obtain for
the integral in the left-hand side of \eqref{reverse_ineq}
\begin{align}\label{lhs}
& \int_{\R} \frac{f_k(y)^2}{1 + \frac{\alpha}{2}y^2}\,dy \ =\
\,\int_{\R}\frac{\chi_k(y)^2}{\left(1 +
\frac{\alpha}{2}y^2\right)^2}\,dy
\\
&  \longrightarrow \ \int_{\R} \frac{1}{\left(1 +
\frac{\alpha}{2}y^2\right)^2}\ dy\ =\ \frac{\pi}{2}
\sqrt{\frac{2}{\alpha}}\ \quad\text{as}\quad k \to \infty\  .
\notag
\end{align}

In view of \eqref{rhs2}, \eqref{lhs},  we obtain
\[
\int_{\R} \frac{f_k(y)^2}{1 + \frac{\alpha}{2}y^2}\,dy \ - \
\frac{1}{2\alpha} \int_{\R} \left(1 + \frac{\alpha}{2}y^2\right)
f_k'(y)^2\,dy \ \longrightarrow\ \frac{3\pi}{8}\
\sqrt{\frac{2}{\alpha}}\ .
\]

From this the conclusion readily follows.

\end{proof}

\medskip

As an immediate consequence of Lemma \ref{L:instab} and Fubini's
theorem we obtain.

\medskip

\begin{cor}\label{C:instab}
For $k\in \mathbb N$, $k\geq k_0$, define $u_k(y,t) = f_k(y)
\chi_k(t)$, where $\chi_k$ and $f_k$ and $k_0$ are as in Lemma
\ref{L:instab}. One has
\begin{equation}\label{E:instab}
\int_{\R^2} \frac{u_k(y,t)^2}{(1 + \frac{\alpha}{2}y^2)(1 +
(\alpha t + \beta)^2)^{3/2}}\,dydy \ > \ \frac{1}{2\alpha}
\int_{\R^2} \frac{1 + \frac{\alpha}{2}y^2}{(1 + (\alpha t +
\beta)^2)^{3/2}} \left(\frac{\partial u_k(y,t)}{\partial
y}\right)^2\,dydt \ .
\end{equation}
\end{cor}

\medskip

We are finally ready to give the

\begin{proof}[\textbf{Proof of Theorem \ref{T:permin}}]
Let $u_{k}$, $k\geq k_0$, be as in Corollary \ref{C:instab}.
Define $a_k:\HH\to\R$, $a_k\in C^\infty_0(\HH)$,  as follows
\[
a_k(x,y,t)\ =\ \frac{\chi_{k}(y)\chi_{k}(t)\chi_{k}(x -
y\,(\alpha\,t + \beta))}{\sqrt{1 + \frac{\alpha}{2}\ y^2}}\ .
\]

We observe that $a_k(\theta(y,t)) = u_k(y,t) \chi_k(0) =
u_k(y,t)$. At this point, appealing to \eqref{he} in Lemma
\ref{L:hardyequiv} and to Corollary \ref{C:instab}, we conclude
that for every fixed $k\geq k_0$, we have for the deformation $\mS
\to \mS^\la = \mS + \la a_k X_1$
\[
 \mathcal V_{II}^H(\mS;a_k X_1)\ =\ \frac{d^2}{d\lambda^2}P_H(\mathcal S^\lambda)\Bigl|_{\lambda =
0}\ <\ 0\ .
\]

This proves that $\mS$ cannot be a local minimizer of the
$H$-perimeter for compactly supported deformations along $X_1$. In a
similar way, using \eqref{he} in Lemma \ref{L:hardyequiv} and
Corollary \ref{C:instab}, we see that $\mS$ cannot be a local
minimizer for deformations along the horizontal normal $\nuX$. In
particular, since every global minimizer is also a local one, $\mS$
cannot be a global minimizer either.

\end{proof}

\end{document}